\documentclass{article}
\usepackage{amsmath}
\usepackage[latin1]{inputenc}

\font\tengoth=eufm10
\font\sevengoth=eufm7
\font\fivegoth=eufm5
\newfam\gothfam
\textfont\gothfam=\tengoth \scriptfont\gothfam=\sevengoth
\scriptscriptfont\gothfam=\fivegoth

\def\blacksquare{\hbox to .60em{\vrule width .60em height .60em}}

\font\bb=msbm10 
\catcode`\@=12  
\def\na{\nabla}

\def\hb {\hfil \break}
\def\n {\vskip 0.2cm \noindent }
\def\scirc{\,{\raise 0.8pt\hbox{$\scriptstyle\circ$}}\,}
\def\ins{\,{\raise 0.2cm \hbox{ $\scriptstyle \circ$}}\,}

\def\a{\alpha}

\def\l{\lambda}
\def\L{\Lambda}
\def\s{\sigma}

\def\h{\eta}
\def\f{\varphi}

\date{}

\begin{document}
	
\centerline{\bf \large  Courbure des tissus en courbes}

\medskip

\centerline{  Daniel Lehmann}

\bigskip

\n {\bf Abstract }

{\it For any $d$-web by curves in an ambiant $n$-dimensional manifold $(d>n)$, we define a vector bundle and a connection on it whose curvature vanishes iff the web has a maximal $(n-1)$-rank.  We wright  a programm on Maple for computing this curvature in the case $d=n+3$ $($but it could be done in the general case$)$ ; as an example,we recover that the exceptional 6-web $W_{0,6}$ in dimension 3 and his subwebs have   maximal rank.}    

 \section{Introduction  :}
 
 Le contexte étant holomorphe ou  analytique réel, 
 on se donne un $d$-tissu en courbes $\cal W$ dans une variété $\cal U$ de dimension $n$, avec $d>n$, localement défini par $d$ feuilletages  ${\cal F}_\lambda$ en courbes. On notera respectivemnt $T_\lambda$ (resp. $N_\lambda$)  le fibré tangent  (resp; normal) à ${\cal F}_\lambda$ :
 $$o\to T_\lambda\to T({\cal U}) \to N_\lambda\to 0.$$
 [Si $V_\lambda$ désigne un champ de vecteurs partout non nul engendrant $T_\lambda$, on supposera que $n$ quelconques de ces  $d$ champs de vecteurs sont toujours linéairement indépendants en chaque point de $\cal U$]. \hb
 La suite $$o\to N^*_\lambda\to T^*({\cal U}) \to T^*_\lambda\to 0$$
 obtenue par dualité permet d'identifier $N^*_\lambda$ aux 1-formes ${\cal F}_\lambda$-semi-basiques (c'est à dire annulées par tout produit intérieur $\iota_v$ par un vecteur $v$ tangent à ${\cal F}_\lambda$), et identifie de m\^eme $\bigwedge^{n-1} N^*_\lambda$ aux $(n-1)$-formes ${\cal F}_\lambda$-semi-basiques. L'application $(\eta_\lambda)_\lambda\to \sum_\lambda \eta_\lambda$ de $\bigoplus_\lambda \bigwedge_\lambda^{n-1}N^*_\lambda$ dans $\Lambda^{n-1}T^*({\cal U})$ est de rang maximum $n$.
 Le  noyau  $E$ de cette application  est donc un fibré vectoriel de rang $d-n$ : 
 $$0\to  E\to \bigoplus_\lambda \bigwedge^{n-1}N^*_\lambda\to \bigwedge^{n-1}T^*({\cal U}).$$
 
 Rappelons qu'une $p$-forme $\omega$ sur $\cal U$ est dite \emph{${\cal F}$-basique} relativement à un feuilletage ${\cal F}$ sur $\cal U$,  si 
 elle est à la fois  ${\cal F}$-semi-basique   et
 ${\cal F}$-invariante ($\iota_v\omega=0$ et $L_v\omega=0$ quel que soit le champ de vecteurs $v$ tangent à ${\cal F}$, $\iota_v$ et $L_v$ désignant respectivement le produit intérieur et la dérivée de Lie). 
 
 Rappelons aussi qu'une \emph{$p$-relation abélienne} du tissu  est une famille de $p$-formes $(\eta_\lambda)_\lambda$ telle que\hb  $\sum_\lambda\eta_\lambda=0$, chaque forme $\eta_\lambda$ étant supposée \emph{${\cal F}_\lambda$-basique}.
 L'ensemble des $p$-relations abéliennes possède  une structure naturelle d'espace vectoriel. 
 [En outre, l'application 
 $(\eta_\lambda)_\lambda\to (d\eta_\lambda)_\lambda$ est une différentielle sur l'espace vectoriel gradué $Ab^*(\cal W)$ des relations abéliennes, dont on note $H^*_{Ab}(\cal W)$ la cohomologie ($d$ désignant la différentielle extérieure des formes différentielles), mais peu importe dans cet article.]
 
 Pour $p=n-1$, on  dira en abrégé \emph{relations abéliennes} au lieu de   $(n-1)$ relations abéliennes : ce sont donc les sections $\sigma=(\eta_\lambda)_\lambda$ de $E$ qui sont solution de l'opérateur différentiel linéaire du premier ordre
 $D \sigma=0$, où $$D\bigl((\eta_\lambda)_\lambda\bigr)=\bigl(L_{V_\lambda}\eta_\lambda\bigr)_\lambda.$$
 
 Nous inspirant de la méthode  proposée par Hénaut dans le cas des tissus planaires ([H1]),  on se propose de montrer :
 
 - que l'espace $R_{h_0-1}$ des relations abéliennes formelles à l'ordre $h_0-1$ possède une structure naturelle de  fibré vectoriel $\cal E$ de rang  $\sum_{h=0}^{h_0-1}(h_0-h)\begin{pmatrix}n-2+h\\h\end{pmatrix}$, où l'on a posé $h_0=d-n$, 
 
- que les relations abéliennes sont les sections $\s$ de $E$ telles que $\na (j^{h_0-1}\s)=0$,    $\na$ désignant  la dérivation covariante d'une certaine connexion tautologique sur $\cal E$, et leur germe en tout point est entièrement déterminé par leur   $(h_0-1)$-jet en ce point.

On retrouve en particulier  un résultat de Damiano ([D1]), selon lequel le rang maximum du tissu est égal à 
$$\sum_{h=0}^{h_0-1}(h_0-h)\begin{pmatrix}n-2+h\\h\end{pmatrix} \hskip 1cm	\Biggl(=h_0.\begin{pmatrix}d\\h_0\end{pmatrix}-(d-1).\begin{pmatrix}d-1\\h_0-1\end{pmatrix}\Biggr). $$
\pagebreak
Ce rang  maximum étant obtenu  ssi la courbure de la connexion précédente est nulle, on redémontre à titre d'exemple,  qu'il en est bien ainsi pour le $6$-tissu exceptionnel $W_{0,6}$ en dimension 3 (ainsi que pour ses 4 et 5-sous-tissus). Bien que le programme Maple   présenté en appendice soit aisément généralisable à tous les cas tels que $d=n+3$, nous ne l'avons fait tourner que pour $n\leq 3$, les temps de calcul devenant très longs pour $n$ plus grand.   Les tissus $W_{0,n+3}$ (redéfinis ci-dessous en termes de champs de vecteurs) ont été particulièrement étudiés par Damiano  ([D1],[D2]) et  Pirio ([Pi3]). Le fait que le rang de ces tissus soit  maximum quel que soit $n$ a été démontré par  Damiano   si $n$ est pair, et par Pirio   si $n$ est impair (ce résultat avait déjà été annoncé    par Damiano lorsque $n$ est  impair, mais avec une erreur dans la démonstration, que Pirio a  corrigée).    Contrairement à  leurs travaux,  notre méthode ne permet  malheureusement   pas  de distinguer les relations abéliennes dites ''combinatoires'' (ce sont celles engendrées par les relations abéliennes des $(n+1)$
 sous tissus), et la relation dite d'Euler\footnote {et qui, si $n$ est impair,  est en fait une relation combinatoire, d'après Pirio.} (qui généralise la relation du dilogarithme pour le  tissu de Bol ([Bo]) dans le cas $n=2$). 
 
Une telle définition de la courbure des tissus, dont la nullité équivaut à la maximalité du rang, remonte à Blaschke ([BB]) dans le cas  $n=2,d=3$. Différentes généralisations ont été abordées depuis, d'abord pour les tissus planaires  ([H1],[Pa], [Pi1]), ensuite  pour les tissus de codimension un ([CL],[DL1], [DL2],[DL3])   lorsqu'ils sont    ''ordinaires'', et aussi pour les tissus en courbes lorsque $d=n+1$ ([DL4]),  ainsi qu'en codimension arbitraire sous des hypothèses assez restrictives ([H2],[L]).

\section{Calcul de $Ab^{n-1}(\cal W)$ :}

On supposera le tissu localement défini par  $d$ champs de vecteurs $V_\lambda$  ($1\leq \lambda\leq d$), et l'on utilisera la base $(V_i)_i$ du module des champs de vecteurs formée par les $n$ premiers.

De façon générale, on notera 
$\lambda ,\mu,\nu,...$ (resp. $i,j,k,...$, resp.$a,b,c,...$) des indices variant de 1 à $d$ (resp. de 1 à $n$, resp. de $n+1$ à $d$),   
\n On posera :
$$V_a=\sum_i f_{i,a}\ V_i,$$ 
et $f_{i,j}=0$ ou $1$ selon que $i$ et $j$ sont  distincts  ou égaux.

\n {\bf Remarque :} {Les composantes $f_{i,a}$ ne s'annulent  jamais  (si, par exemple,  $f_{1,a}$ était nul, les champs ($V_a$, $V_2$,...,$V_n$) seraient linéairement dépendants). 
	
\n On peut donc supposer (ce que nous ferons désormais) :
$$f_{n,a}\equiv 1 \hbox{ pour tout } a.  $$	
\n On notera :
 
 $(\a_i)$ la base du module des 1-formes locales duale de la base  $(V_i)$ du module des champs de vecteurs locaux,
 
 $\Lambda$ la $n$-forme volume associée $\Lambda:= \a_1\wedge \a_2\wedge ...\wedge\a_{n-1}\wedge \a_n $,
 
  $\L_i$ la $(n-1)$-forme $\Lambda_i:=(-1)^{i+1} \a_1 ...\wedge \widehat{ \a_i}\wedge ...\wedge \a_n $ obtenue (au signe près) en omettant le terme $\a_i$:$$\Lambda_i:=(-1)^{i+1} \a_1 ...\wedge \widehat{ \a_i}\wedge ...\wedge \a_n. $$
  
  $\Lambda_{ij}$ (avec $i<j$) la $(n-2)$-forme obtenue (au signe près) en omettant les termes $\a_i$ et $\a_j$ $$\Lambda_{ij}:=(-1)^{i+j}\ \a_1 ...\wedge \widehat{ \a_i}\wedge... \wedge \widehat{ \a_j}...\wedge \a_n. $$
 
 \n Notant $C_{ij}^k$ les fonctions telles que $[V_i,V_j] =\sum_k C_{ij}^k V_k$, on obtient les formules :

$d\a_k=-\sum_{i<j} C_{ij}^k\ \a_i\wedge \a_j$,

   et  $d\L_i=\varphi_i\  \L$ avec 
  $\varphi_i:=\sum_{j,j<i}C_{ji}^j-\sum_{j,i<j}C_{ij}^j$.

\n {\bf Lemme 1:}

{\it  Une famille $(\eta_\lambda)_\lambda$ de $(n-1)$ formes différentielles
$ \eta_\lambda=\sum_i X_{i,\lambda}\Lambda_i$ est une $(n-1)$ relation abélienne, ssi sont vérifiées les équations suivantes :

	$(I)$  \hskip 1cm    $X_{i,j}=0$ si $i\neq j$,\hskip 1cm  et  \hskip 1cm  $X_{i,a}=f_{i,a}X_{n,a}$ pour tout $(i,a)$,
	
	$(II )$ \hskip 1cm   $ X_{i,i} =-\sum_a X_{i,a}=0  $  pour  tout $i$, 
	
	$(III )$  \hskip 1cm $\sum_a V_i.(f_{i,a}X_{n,a})=-\f_i.\sum_a ( f_{i,a}X_{n,a}) $  pour  tout $i$,
	
\hskip 1.5cm	 et \ $\sum_i V_i.(f_{i,a}.X_{n,a})=-\Bigl(\sum_i \f_i f_{i,a}\Bigr).X_{n,a}. $  pour  tout $a$.}

\n {\it Démonstration :} 

\n Dire que les produits intérieurs $\iota_{V_\lambda}(\eta _\lambda)$ sont tous nuls s'écrit en effet :
 $$ (I_{ij,\lambda})\hskip 1cm f_{i,\lambda}X_{j,\lambda}-f_{j,\lambda}X_{i,\lambda}=0 \hbox{ pour tout couple $(i<j)$ et tout $\lambda$}.$$
 En particulier, $(I_{ij,j})$ s'écrit $X_{i,j}=0$ si $i\neq j$, et $(I_{in,\a})$ s'écrit $X_{i,a}=f_{i,a}X_{n,a}$. Toutes les autres équations $(I_{ij,\lambda})$ en découlent. 
 
\n Les équations $(II)$ expriment que la somme des formes $\eta_\lambda$ est nulle.

\n Les équations $(III)$ expriment, compte tenu de $(I)$ et $(II)$, que les dérivées de Lie $L_{V_\lambda}( \eta _\lambda)$ sont toutes nulles.

\rightline{QED}

\n{\bf Inconnues principales :} 

Finalement, les équations  $(I)$ et $(II)$ montrent que les $n.d$ fonctions inconnues  $X_{i,\l}$ ne sont pas indépendantes et que l'on peut   peut finalement se ramener aux $d-n$ inconnues  principales indépendantes 
$Y_a:=X_{n,a}$ : ce  sont les composantes  d'une section $\sigma$ du  fibré vectoriel $E$, relativement à une trivialisation ad-hoc de ce fibré.

\n{\bf L'opérateur différentiel :} 

Sommant terme à terme les équations $(III)$, on obtient : $\sum_{\lambda}\sum_i \bigl(X_{i,\lambda}\bigr)'_i=0$, soit
$\sum_{i}\bigl(\sum_\lambda X_{i,\lambda}\bigr)'_i=0$. Puisque $\sum_\lambda X_{i,\lambda}=0$, d'après $(II)$, on en déduit que 
 les $d$ équations de $(III)$ ne sont pas linéairement indépendantes : 
 l'une d'elles (disons la dernière pour fixer les idées) équivaut à  la somme des $d-1$ autres, et peut donc \^etre omise.

En résumé, $(III)$ permet donc d'exprimer l'ensemble  $R_1$  des  relations abéliennes formelles à l'ordre 1 comme les éléments  de $J^1E$ annulés par un opérateur différentiel $D_1$ linéaire d'ordre 1 :
$$ 0\to R_1\to  J^1E\buildrel{D_1}\over\longrightarrow \Bigl(\wedge^{n-1}T^*({\cal U})\Bigr)^{\oplus d}  $$
Réécrivant $(III)$ en isolant le symbole principal $P_1$ de $D_1$, les éléments de $R_1$    se projetant sur un  0-jet donné dans $E$  sont  solutions d'un système linéaire $\Sigma_1$ de $d-1$ équations à $n\times (d-n)$ inconnues $(X_{n,a})'_i$, avec second membre 
 
 $(III_i)$ \hskip 1.5cm $\sum_a   f_{i,a}.(Y_a)'_i=-\sum_a \Bigl((f_{i,a})'_i+\f_i.f_{i,a}\Bigr)Y_a $  pour  tout $i$,

 $(III_a)$\hskip 1cm	 et \ $\sum_i   f_{i,a}.(Y_a)'_i=-\Bigl(\sum_i \bigl((f_{i,a})'_i+\f_i.f_{i,a}\bigr)\Bigr)Y_a $  pour   $n+1\leq a\leq d-1$.}

\n $(...)'_i$ désignant la dérivée d'une fonction par rapport au champ de vecteurs $V_i$.
\n 
On écrit $\Sigma_1$  sous la forme :
$$ <P_1\ ,\ Y^{(1)}>+<Q_1\ ,\ Y>=0 \hbox{ ou }  <M_1\ ,\ j^1Y>=0 \hbox{ en abrégé,  } $$
où l'on   note :

$Y$ le $(d-n)$ vecteur de composantes $\bigl((Y_a)\bigr)_{a}$\ ,

  $Y^{(1)}$ le $(d-n)\times n$ vecteur de composantes $\bigl((Y_a)'_i\bigr)_{(a,i)}$\ ,

$P_1$  la matrice $((f_{i,a}))$ de taille $(d-1)\times n(d-n$), correspondant à la partie homogène  des équations $(III_\l)$\ ,

et $Q_1$   la matrice de taille $(d-1)\times (d-n)$ correspondant aux  seconds membres de ces  équations.

 \n {\bf Lemme 2 :} {\it La matrices $P_1$ est  de rang maximum $d-1$ en tout point de $\cal U$, et $R_1$ est donc un fibré vectoriel de rang $(d-n)+\bigl(n.(d-n)-(d-1)\bigr)$.
}

\n {\it Démonstration :}

On conviendra  d'ordonner les $(d-n)\times n$ composantes $(Y_a)'_i$ de $Y^{(1)}$
en écrivant d'abord les composantes $(Y_a)'_1$ avec $n+1\leq a\leq d$, puis les composantes $(Y_a)'_2$,..., puis les composantes $(Y_a)'_n$. 
La matrice $P_1$, de taille $(d-n).n\times (d-1)$ se décompose alors  en $n$ blocs $T_i$ de taille $(d-n)\times (d-1)$
corespondant aux $n$ dérivées premières du vecteur $Y$, de sorte que $P_1=\left(\begin{array}{ccc}T_1&...&T_n\end{array}\right)  $  et $$<P_1,Y^{(1)}>=\sum_{i=1}^n <T_i,Y'_i>.$$ 
On décomposera aussi  chaque bloc $T_i=\left(\begin{array}{c}A_i\\B_i\end{array}\right)$    en 2 blocs superposés $A_{i}$ et $B_i$,  le premier (de taille \hb $(d-n)\times (n-1)$)  correspondant aux équations $III_\l$ pour 
$1\leq \l\leq n-1$, et le second (de taille $(d-n)\times (d-n)$) correspondant aux équations $III_\l$ pour 
$n\leq \l\leq d-1$ en convenant de les écrire dans l'ordre $(n+1,n+2,...,d-1,n)$ avec $n$ en dernier.
On redécompose enfin $P_1$ en 4 blocs de la façon suivante :
$$P_1=\left(\begin{array}{cc}
	\bigl(A_1\ ....\ A_{n-1}\bigr)	& \bigl(­A_n =0\bigr)\\
	\bigl(B_1\ ....\ B_{n-1}\bigr)	& \bigl(B_n\bigr)
\end{array}\right).$$
La	$i$-ième ligne du  bloc $A_i$ ($1\leq i\leq n-1$) s'écrit $(f_{i,n+1} ,f_{i,n+2},...,f_{i,d})$, et toutes les autres lignes n'ont que des $0$. En particulier, le bloc $A_n$ n'a que des zéros. Quand à la matrice carrée $B_n$, elle  n'a que des 1 sur la dernière ligne, et que des 0 sur les autres lignes sauf à la $i$-ème colonne où le coefficient est 1. 

\rightline{QED}

\n [En résumé, indexant par $(i,a)$ les $n.(d-n)$ colonnes de $P_1$ et par $\l$, ses $(d-1)$ lignes $(1\leq\l\leq d-1)$, le coefficient $(P)_\l^{i,a}$ de $P_1$ sur la ligne $\l$ et la colonne $(i,a)$ est égal à $f_{i,a}$ lorsque $\l=i$ ou $\l=n+i$, et à 0 sinon, avec $f_{n,a}=1$].
  
\rightline{QED}

\n Exemples pour $n=3$, $d=6$ :

$$P_1=\left(\begin{array}{ccccccccccc}
	f_{1,4}	&f_{1,5}  & f_{1,6}& & 0 & 0 & 0 && 0 & 0 & 0 \\
	0	& 0 & 0 & &f_{2,4}	&f_{2,5}  & f_{2,6} & & 0 & 0 & 0  \\
	&&&&&&&&&&\\

	f_{1,4}	& 0 & 0 &&  f_{2,4}& 0 &0  &&1  & 0 & 0  \\
	0	& f_{1,5} &0 & & 0 & f_{2,5} &0  && 0 & 1&  0\\  
	0 & 0 & 0 &&0  & 0 & 0 && 1	&1  &1 \\ 
\end{array}\right)$$


\n  Puisqu'aucun des coefficients $f_{i,a}$ n'est nul,  le bloc $\bigl(A_1\ ....\ A_{n-1}\bigr)$ de $P_1$ est de rang $n-1$, tandis que $B_n$ est de rang $d-n$. Le bloc $A_n$ n'ayant que des zéros, le rang de $P_1$ est égal à la somme $(n-1)+(d-n)=d-1$ des rangs 
des blocs diagonaux. 

\n {\bf Les prolongements de l'opérateur différentiel  :} 

  On notera  $R_h$ le sous-ensemble   des $h$-jets de $E$
constitué  des 
relations abéliennes formelles à l'ordre $h$ : c'est l'ensemble des solutions du $(h-1)$-ème  prolongement $D_h$ de l'opérateur différentiel $D_1$ du premier ordre  défini ci-dessus, et se définit par récurrence par la formule 
$$R_h= J^hE\cap J^1R_{h-1} \hskip 1cm \hbox{$\bigl($intersection dans $J^1(J^{h-1}E)\bigr)$},$$
qui prendra un sens dès que l'on aura démontré que $R_h$ possède une structure de  sous-fibré vectoriel de $J^hE$, ce qui est un corollaire du lemme 4 ci-dessous. 

Soit $I=(i_1,i_2,...,i_n)$ un multi-indice  de dérivation à l'ordre $h$ (avec $i_1+i_2+...+i_n=h$). On conviendra que la dérivée  $f'_I$ d'ordre $h$ d'un fonction $f$ signifiera que l'on a d'abord dérivé la fonction $j_1$ fois par rapport à $V_1$, puis $j_2$ fois par rapport à $V_2$,..., puis $j_n$ fois par rapport à $V_n$. 

\n {\bf Remarque : } Puisque  $V_i.(V_j.f)=V_j.(V_i.f)+[V_i,V_j].f$,   toute dérivation d'ordre $h$ effectuée dans un ordre différent peut \^etre remplacée par la dérivation précédente,  modulo des dérivations d'ordre strictement inférieur, ce qui ne modifie  pas le symbole principal de l'opérateur différentiel $D_h$ d'ordre $h$.

Notant alors  $$c(n,h)=\frac{(n-1+h)!}{(n-1)!\ h!} $$  la dimension de l'espace vectoriel des polyn\^omes homogènes de degré $h$ à $n$ variables, l'espace affine des   relations abéliennes formelles à l'ordre $h$  se  projetant sur une  relation abélienne formelle à l'ordre $h-1$ donnée sont   les solutions   d'un système linéaire $\Sigma_h$ de $c(n,h-1).(d-1) $ équations à $  c(n,h).(d-n)$ inconnues $((Y)'_I)_{|I|= h}$, avec second membre, écrit sous la forme :
$$ <P_h\ ,\ Y^{(h)}>+<Q_h\ ,\ j^{(h-i)}Y>=0 \ , $$
 où 

$Y^{(h)}$ désigne le $ (d-n)\times c(n,h)$ vecteur de composantes  $\bigl((Y_a)'_I)_{(a,i)}$, 

$P_h$ désigne  l'opérateur linéaire  de taille $(d-1).c(n,h-1)\times (d-n).c(n,h)$ représentant le symbole principal de $D_h$, 

 $Q_h$   l'opérateur linéaire de taille $(d-1)\times c(n,h-1)$ correspondant au   second  membre du système. 
 
\n Notant $j^h Y$ le $(d-n)$ vecteur $\oplus_{k=0}^h Y^{(k)}$, il pourra  aussi être utile de regrouper tous les systèmes linéaires $\Sigma_k$ ($k\leq h$) sous la forme d'un système linéaire homogène unique $$<M_h,j^h Y>=0$$ à
$c(n+1,h-1).(d-1) $ équations à $  c(n+1,h).(d-n)$ inconnues $((Y)'_I)_{|I|\leq h}$, les opérateurs  
$M_h=\Bigl(\!\Bigl(M_{H,u}^{K,v}\Bigr)\!\Bigr)$, emboités les uns dans les autres, sont tels que l'équation $(III_u)'_H$ s'écrit :
$$\sum_{K,v}M_{H,u}^{K,v}.(Y_{n+v})'_K=0 \hbox{ avec } |K|\leq h \hbox{ et } 1\leq v\leq d-n.$$

 \n Les coefficients $M_{H,u}^{K,v}$ de $M_h$ (donc aussi ceux de $P_h$ et $Q_h$) s'obtiennent par récurrence sur $|H|$ par dérivations successives de la formule  
$ <P_1\ ,\ Y^{(1)}>+<Q_1\ ,\ Y>=0$.
$$M_h=\begin{pmatrix}M_{h-1}&0\\
	Q_{h}&P_{h}\
\end{pmatrix} $$

Précisons les trivialisations permettant d'exprimer tous les opérateurs précédents en termes de matrices. A cet effet, on  commence par numéroter les $n$-multi-indices de dérivation. 
Pour un tel multi-indice   $I=(i_1,i_2,...,i_n)$, on notera   $|I|$  son ordre $\sum_{k=1}^n i_k$. On numérote tous ceux de ces indices qui sont d'ordre au plus $h_0=d-n$ à l'aide d'une bijection $LL:\{0,1,\cdots ,c(n+1,h_0)-1\}\to \{multi-indices\}$ de façon que 

$LL(0)=(0,0,...,0)$, 

$LL(k)$ désigne, pour tout $k$  ($1\leq k\leq n$),  le multi-indice $1_k=(0,0,...,0,1,0,...,0)$ avec 1 à la k-ième place,

 et 
que les $c(n,h)$ multi-indices  d'ordre $h$ soient plus généralement numérotés entre $c(n+1,h-1)+1$ et $c(n+1,h)$.

 \n On indexera  par $i=u+(d-1).t$ la ligne  de $M_h$  correspondant à la dérivée $(III_u)'_H$ d'indice $H=LL(t)$  de l'équation $III_u$\ ,    avec $|H|\leq h-1$ et  $1\leq u\leq d-1$. On écrira aussi $ (III_u)'_t$
 
\n  De m\^eme, on indexera par  $j=v+(d-n).s$ la colonne   de $M_h$  correspondant à la dérivée $(Y_{n+v})'_K$ d'indice $K=LL(s)$ (notée aussi $(Y_{n+v})'_s$),   avec $|K|\leq h$ et  $1\leq v\leq d-n$.

\n Réciproquement, les fonctions $i\to H(i) =LL([\frac{i+d-2}{d-1}])$  et  $i\to u(i) =i-[\frac{i+d-2}{d-1}]$ permettent de retrouver $(III_u)'_H$ à partir de $i$, 
de m\^eme que les fonctions $j\to K(j) =LL([\frac{j+d-n-1}{d-n}])$  et  $j\to v(j) =j-[\frac{j+d-n-1}{d-n}]$ permettent de retrouver $(Y_{n+v})'_K$ à partir de $j$, et l'on définit la matrices représentant $M_h$ par les coefficients 
$$F(i,j) =M_{H(i),u(i)}^{K(j),v(j)}$$  ou $$ F\bigl(u+(d-1)t,v+(d-n)s)\bigr) =M_{H,u}^{K,v} \hbox{\hskip .5cm avec \hskip .5cm $H=LL(t)$ et $K=LL(s)$}.$$ 
Ces coefficients se calculent par récurrence sur l'ordre $h\bigl(LL(t)\bigr)$ (noté aussi $h(t)$ plus simplement) : supposant par exemple, pour simplifier, que  $V_i=\frac{\partial }{\partial x_i}$ pour $i$ compris entre 1 et $n$ relativement à un certain système $(x_i)$ de coordonnées locales : les équations initiales $E_u$  s'écrivent  $$\sum_v\Bigl(F(u,v).Y_v+\sum_{m=1}^n F\bigl(u,v+(d-n).m\bigr).(Y_v)'_m\Bigr)\equiv 0.$$
Les coefficients 
$F(u,v)$ et $F(u,v+(d-n).m)$ correspondant au cas $t=0$ sont connus. Plus généralement, supposons déjà connus les coefficients
$F(u+(d-1)t, v+(d-n)s) $ de l'équation 
$$(III_u)'_t\hskip 2cm \sum_v\Bigl(F(u+(d-1)t,v).Y_v+\sum_{s,1\leq h(s)\leq h(t)+1}F\bigl(u+(d-1)t,v+h_0.s\bigr).(Y_v)'_s\Bigr)\equiv 0.$$
Pour tout entier $k$ ($ 1\leq i\leq n$), on notera  $ad_k(t)$ le multi-indice tel que $LL\bigl(ad_k(t)\bigr)$ s'obtient à partir du multi-indice  $LL(t)$ en augmentant de une unité son $k$-ième terme, sans toucher aux autres : l'équation 
$(III_u)'_{ad_k(t)}$  obtenue en dérivant $(III_u)'_t$
 par rapport à $x_k$ s'écrit alors :
$$ 
\sum_v\biggl((F(u+(d-1)t,v))'_k.Y_v+\sum_{s,h(s)\geq 1}\Bigl(F\bigl(u+(d-1)t,v+h_0.s\bigr)\Bigr)'_k.(Y_v)'_s+
F\bigl(u+(d-1)t,v+h_0.s\bigr).(Y_v)'_{ad_k(s)}\biggr)\equiv 0,
$$d'où les formules de récurrence :

$F\bigl(u+(d-1).ad_k(t),v\bigr)=\Bigl(F\bigl(u+(d-1)t,v\bigr)\Bigr)'_k$,

$F\bigl(u+(d-1)ad_k(t),v+h_0.s\bigr)=\Bigl(F\bigl(u+(d-1)t,v+h_0.s\bigr)\Bigr)'_k$ si $h(s)\geq 1$ et si $s$ n'appartient pas à l'image de l'application $ad_k$,

$F\bigl(u+(d-1).ad_k(t),v+h_0.ad_k(r)\bigr)=\Bigl(F\bigl(u+(d-1)t,v+h_0.ad_k(r)\bigr)\Bigr)'_k+F\bigl(u+(d-1)t,v+h_0.r\bigr)$  si $0\leq h(r)\leq h(t)$,

$F\bigl(u+(d-1).ad_k(t),v+h_0.ad_k(r)\bigr)=F\bigl(u+(d-1)t,v+h_0.r\bigr)$  si  $ h(r)= h(t)$,

$F\bigl(u+(d-1).ad_k(t),v+h_0.s\bigr)=0$ si $ h(s)= h(t)+1$ et  si $s\notin Im(ad_k)$ ou   si $ h(s)> h(t)+1$

\n {\bf Lemme 3 :} 

{\it Posant $h_0=d-n$, le système linéaire $\Sigma_h$ a plus $($resp. moins$)$ d'inconnues que d'équations selon que $h<h_0$ $($resp. $h>h_0)$. Et il en a autant pour $h=h_0$. }

\n C'est une conséquence immédiate de la formule 
$$h_0 .c(n,h) -(d-1). c(n,h-1)=(h_0-h).c(n-1,h),$$aisée à vérifier. 

\n {\bf Lemme 4 :} {\it Tous les systèmes linéaires  $\Sigma_h$, $($c'est-à-dire toutes les matrices $P_h)$  sont   de rang maximum, à savoir :
	
	$(d-1). c(n,h-1)$  si $h\leq h_0$,
	
	$h_0 .c(n,h)$  si $h\geq h_0$.}

On observe en effet que   $P_h$ se décompose en $c(n,h).c(n,h-1)$ blocs $K_J^I$ de taille $h_0.(d-1)$ où $J=(j_1,...,j_n)$ désigne un multi-indice d'ordre $|J|=h-1$, et $I=(i_1,...,i_n)$ désigne un multi-indice d'ordre $|I|=h$, $K_J^{I'}$ étant plus  à droite que $K_J^{I}$ (pas nécessairement juste à c\^oté) si $I'>I$, et $K_{J'}^{I}$ étant en dessous (pas nécessairement juste en dessous) de $K_J^{I}$ si $J'>J$. De façon précise, 

$K_J^I=T_k$ si $i_k=j_{k}+1$ et $i_\ell=j_{\ell}$ pour $\ell\neq k$   et l'on écrira alors $I=J+1_k$, 

$K_J^I=0$ n'a que des zéros sinon.

\n Exemple pour $n=3$ et $d=6$  : \hskip 1cm $ P_1=\left(\begin{array}{ccc}
	A_1& A_{2}	& \\
	B_1 & B_{2}	& B_3
\end{array}\right),$

 $$ \ P_2=	\left(\begin{array}{cccccc}
A_1&A_2&&&&\\ 
B_1&B_2&B_3&&&\\
&A_1&&A_2&&\\ 
&B_1&&B_2&B_3&\\
& &A_1&&A_2&\\ 
&&B_1&&B_2&B_3\\	
\end{array}\right),    \ P_3=	\left(\begin{array}{cccccccccc}
A_1	& A_2 &  &  &  &  &  &  &  &  \\     
B_1	& B_2 & B_3 &  &  &  &  &  &  &\\     
& A_1 &  & A_2 &  &  &  &  &  &\\     
& B_1 &  & B_2 & B_3 &  &  &  &  &\\ 
&  &A_1  &  & A_2 &  &  &  &  &\\     
&  &B_1  &  &B_2  & B_3 &  &  &  &\\  
&  &  &A_1  &  &  & A_2 &  &  &\\     
&  &  & B_1 &  &  &B_2  & B_3 &  &    \\
&  &  &  & A_1 &  &  & A_2 &  &  \\
&  &  &  & B_1 &  &  & B_2 & B_3 &   \\
&  &  &  &  &A_1  &  &  &A_2  &      \\
&  &  &  &  &B_1  &  &  &B_2  &B_3  
\end{array}\right).$$
[Il est sous-entendu que les blocs laissés en blanc sont des blocs de zéros]. 

 On remarque que, pour $h\leq h_0$,  l'indice colonne  maximum des  composantes non nulles du   vecteur  ligne $(J,\l)$   (où $J=(j_1,...,j_n)$ désigne  un multi-indice d'ordre  $ h-1$), est une fonction strictement croissante de $(J,\l)$ quand on se limite aux indices lignes tels que $\l\neq 1$ : ces vecteurs lignes sont donc nécessairement 
linéairement indépendants. Compte tenu de leurs autres composantes non nulles, on voit que les vecteurs lignes $(J,1)$ sont linéairement indépendants des autres. 

De même, si $h\geq h_0$,  l'indice ligne  maximum des  composantes non nulles du   vecteur  colonne $(I,\l)$   (où $I=(i_1,...,i_n)$ désigne  un multi-indice d'ordre  $ h$), est une fonction strictement croissante de $(I,\l)$ quand on se limite aux indices colonnes  tels que $\l\neq d-1$ : ces vecteurs colonnes sont donc nécessairement 
linéairement indépendants. Compte tenu de leurs autres composantes non nulles, on voit que les vecteurs colonnes $(I,d-1)$ sont linéairement indépendants des autres. 

   \rightline{QED}

\n {\bf Théorème 1 :}

{\it 
	
		$(i)$ Posant $h_0=d-n$,  $R_h$ est un fibré vectoriel de rang  $\sum_{j=0}^{h-1}(h_0-j).c(n-1,j)$ si $h\leq h_0-1$,\hb \indent  et 
	 est au plus égal à $\sum_{j=0}^{h_0-1}(h_0-j).c(n-1,j)$ si $h\geq h_0$.
	 
	$(ii)$ le $(n-1)$-rang du tissu est au plus égal à $$\sum_{j=0}^{h_0-1}(h_0-j).c(n-1,j) 
\ \	\Bigl(=h_0.c(n+1,h_0 )-(d-1).c(n+1,h_0-1)\Bigr). $$
	
		$(iii)$ On définit une connexion sur ${\cal E}:=R_{h_0-1}$ en posant\footnote{L'expression $(s,<..>)$ appartient à   $J^1J^{h_0-1}E$, dans la mesur où un élément de ce fibré est       défini par sa partie principale et   sa projection sur $J^{h_0-1}E$.   D'autre part, la suite exacte ci-dessus ainsi que les systèmes linéaires $(\Sigma_{h_0}, P_{h_0},Q_{h_0})$ ayant  une signification intrinsèque, la connexion ainsi définie ne dépend   pas de l'ordre des champs de vecteurs $V_\l$. }, pour toute section $s$ de $R_{h_0-1}$ :  $$\nabla s=j^1s-\Bigl(s\ , <-(P_{h_0})^{-1}. Q_{h_0},s>\Bigr).$$
		
	$(iv)$	 Les $(n-1)$-relations abéliennes 
		du tissu sont alors les sections $\sigma$ de $E$ telles que  $\nabla (j^{h_0-1}\sigma) =0$.
		
		$(v)$ Le tissu est de rang maximum ssi la courbure de la connexion $\nabla $ est nulle. 	

}

\n {\bf Remarque :} La conclusion  	$(ii)$ est un résultat de Damiano ([D]), que l'on redémontre  ici. 

\n {\it Démonstration :}

\n La matrice $P_h$ étant de rang maximum, 

- les solutions de $\Sigma_h$ forment     un espace affine de dimension  $(d-n-h).c(n-1,h)$ si $h\leq h_0$, 

- le système est cramérien pour $h=h_0$ et est surdéterminé si $h>h_0$. 

\n On en déduit $(i)$. Puisque le contexte est analytique, un germe de relation abélienne est entièrement défini par sont jet infini : $(ii)$ en résulte. 

\n La   matrice $P_{h_0}$ étant inversible, l'application 
$s\to \Bigl(s\ , <-(P_{h_0})^{-1}.Q_{h_0},s>\Bigr)$ définit une scission $\tau$ de la suite exacte 
 $$0\to T^*{\cal U}\otimes R_{h_0-1}\to J^1 R_{h_0-1}\buildrel{\buildrel{\tau}\over{\longleftarrow}}\over{\longrightarrow}  R_{h_0-1}\to 0 $$et, pour toute section $s$ de  $R_{h_0-1}$, $\tau(s)$ et $j^1s$ ont m\^eme projection sur $R_{h_0-1}$  : leur différence est donc la  dérivation covariante d'une connexion sur $R_{h_0-1}$, d'où $(iii)$.
  
  \n Enfin,  dire que $j^1j^{h_0-1}\sigma= \Bigl(j^{h_0-1}\sigma\ ,\ <(P_{h_0})^{-1}.Q_{h_0},j^{h_0-1}\sigma>\Bigr)$ signifie que  le  $(h_0-1)$-ième prolongement de l'opérateur différentiel $D_1$ annulle $\sigma$, d'où 	$(iv)$ et 	$(v)$.
 
\rightline{QED}

\n {\bf Cas particulier :} On suppose que  $V_i=\partial_i$ pour tout $i$, où l'on s'est donné un système de coordonnées locales $(x_i)$, et où l'on a écrit en abrégé $\partial_i$ au lieu de $\frac{\partial }{\partial x_i}$. On a alors $f_{i,j}=0$ ou $1$, selon que $i$ est différent ou égal à $j$. Les fonctions $\varphi_i$ sont identiquement nulles.   Les équations 
$(III_i)$ et $(III_a)$ se simplifient et  deviennent respectivement :
$$\sum_a (f_{i,a} .Y_a)'_i=0 , \hskip 1cm\hbox{ et }\hskip 1cm  \sum_i (f_{i,a} .Y_a)'_i=0 .$$
 {\bf Exemple du $(n+3)$-tissu exceptionnel $W_{0,n+3}$ :}

Relativement à un système $(x_i)$ de coordonnées affines dans l'espace projectif, on définit, pour $x_i\neq 0,1$ : 
\n $V_i=\partial_i$ pour tout $i$, \hskip .4cm 
$V_{n+1}= \frac{1}{x_n}\sum_i x_i\ \partial_i$,\hskip .4cm 
$V_{n+2}=\frac{1}{x_n-1}\sum_i (x_i-1)\ \partial_i$,\hskip .4cm 
$V_{n+3}=\frac{1}{x_n(x_n-1)}\sum_i x_i(x_i-1)\ \partial_i$ :\hskip .4cm 

\n Le $(n+3)$ tissu ainsi défini est bien égal au tissu $W_{0,n+3}$ (notations de [Pi3]), car

- les intégrales premières  $x_j$ $(j\neq i)$ du feuilletage ${\cal F}_i$ sont invariantes par  $V_i$,

- les  intégrales premières $x_j/x_n$ $(j\neq n)$ du feuilletage ${\cal F}_{n+1}$ sont invariantes par $V_{n+1}$,

- les intégrales premières $(x_j-1)/(x_n-1)$ $(j\neq n)$ du feuilletage ${\cal F}_{n+2}$ sont  invariantes par $V_{n+2}$,

- les intégrales premières $x_j(x_n-1)/x_n(x_j-1)$ $(j\neq n)$ du feuilletage ${\cal F}_{n+3}$ sont invariantes par $V_{n+3}$.

\n Une autre présentation de ce tissu en termes de fonctions oubli sur l'espace ${\cal M}_{0,n+3}$ des configurations de $n+3$ points de {\bb P}$^1$ (cf. [Pi3]), bien que moins maniable techniquement dans la situation qui nous intéresse,   présente l'avantage de faire jouer le m\^eme rôle à tous les $(n+3)$ feuilletages du tissu, ce qui implique en particulier que, pour un entier $d$ donné $(n+1\leq d<n+3)$,   tous les $d$-sous-tissus de $W_{0,n+3}$ sont isomorphes.

 Redémontrons, au moins\footnote{En théorie, le programme proposé fonctionne quel que soit $n$, mais les temps de calcul deviennent très longs avec un simple ordinateur portable.} pour $n=3$,   que le rang de ce tissu est maximum en vérifiant, à l'aide d'une programmation sur Maple,  que la  courbure définie ci-dessus est nulle. Nous allons en fait faire le calcul pour un tissu $W_c$ dépendant d'un paramètre scalaire $c$, et égal à $W_{0,6}$ pour $c=0$ : on obtient alors une courbure dont certains  coefficients (notés $O(c)$) sont non nuls, mais  petits d'ordre $c$ quand $c$ tend vers 0. Plus précisément, on obtient à la fin les trois matrices  $10\times 10$ de courbure  dont tous les coefficients sont nuls ou $O(c)$. [Le tissu $W_c$ est simplement obtehu à partir du 6-tissu $W_{0,6}$ en remplaçant  le quatrième feuilletage (le bouquet des droites passant par l'origine) par le bouquet des droites pasant par le point $(-c,-c,-c)$ : les cubiques rationnelles qui sont les feuilles du sixième feuilletage ne passent génériquement plus par ce point si $c\neq 0$, et $W_c$ n'est pas isomorphe à $W_{0,6}$, bien que le  5-sous tissu obtenu en omettant le sixième feuilletage soit isomorphe à n'importe lequel des 5-sous-tissus de $W_{0,6}$].

Ce programme, adapté ici au cas $n=3$, se généralise aisément à tous les tissus tels que $d=n+3$, mais les temps de calcul peuvent  devenir très longs. On obtient à la fin   les trois matrices de courbure

$$ Ko(1,2)=\begin{pmatrix}
	0&0&0&0&0&0&0&0&0&0\\
	0&0&0&0&0&0&0&0&0&0\\
	0&0&0&0&0&0&0&0&0&0\\
	0&0&0&0&0&0&0&0&0&0\\
	0&0&0&0&O(c)&O(c)&O(c)&O(c)&O(c)&O(c)\\
	0&0&0&0&O(c)&O(c)&O(c)&O(c)&O(c)&O(c)\\
	0&0&0&0&O(c)&O(c)&O(c)&O(c)&O(c)&O(c)\\
	0&0&0&0&0&0&0&0&0&0\\
	0&0&0&0&0&0&0&0&0&0\\
	0&0&0&0&0&0&0&0&0&0\\
\end{pmatrix}$$
$$
Ko(1,3)=\begin{pmatrix}
	
	0&0&0&0&0&0&0&0&0&0\\
	0&0&0&0&0&0&0&0&0&0\\
	0&0&0&0&0&0&0&0&0&0\\
	0&0&0&0&0&0&0&0&0&0\\
	0&0&0&0&O(c)&O(c)&O(c)&O(c)&O(c)&O(c)\\
	0&0&0&0&0&0&0&0&0&0\\
	0&0&0&0&0&0&0&0&0&0\\
	0&0&0&0&O(c)&O(c)&O(c)&O(c)&O(c)&O(c)\\
	0&0&0&0&0&0&0&0&0&0\\
	0&0&0&0&O(c)&O(c)&O(c)&O(c)&O(c)&O(c)\\
\end{pmatrix}
$$

$$ Ko(2,3)=\begin{pmatrix}
	0&0&0&0&0&0&0&0&0&0\\
	0&0&0&0&0&0&0&0&0&0\\
	0&0&0&0&0&0&0&0&0&0\\
	0&0&0&0&0&0&0&0&0&0\\
	0&0&0&0&0&0&0&0&0&0\\
	0&0&0&0&0&0&0&0&0&0\\
	0&0&0&0&O(c)&O(c)&O(c)&O(c)&O(c)&O(c)\\
	0&0&0&0&0&0&0&0&0&0\\
	0&0&0&0&O(c)&O(c)&O(c)&O(c)&O(c)&O(c)\\
	0&0&0&0&O(c)&O(c)&O(c)&O(c)&O(c)&O(c)\\
\end{pmatrix}$$

\n {\bf Remarques :} 

Notant $$A=\sum_{k=1}^3 A(k)\ dx_k$$ la forme de la  connexion tautologique définie sur sur ${\cal E} $, et 
 $$ \ Ko=Ko(1,2)\ dx_1\wedge dx_2 +Ko(1,3)\ dx_1\wedge dx_3 +Ko(1,2)\ dx_2\wedge dx_3 $$ sa forme de courbure, on observe que sont  nuls tous les coefficients 

- des   sous-matrices  $6\times 4$    formées avec les 4 premières colonnes et les 6 dernières lignes des matrices de connexion $A(k)$,

- des sous-matrices  $10\times 4$    formées des 4 premières colonnes  des matrices de courbure $Ko(k,m)$.
$$ A(k)=\begin{pmatrix}
	\varphi(k)&w&w&w&w&w&w&w&w&w\\
	0&w&w&w&w&w&w&w&w&w\\
	0&w&w&w&w&w&w&w&w&w\\
	0&w&w&w&w&w&w&w&w&w\\
	&&&&&&&&&\\
	0&0&0&0&w&w&w&w&w&w\\
	0&0&0&0&w&w&w&w&w&w\\
	0&0&0&0&w&w&w&w&w&w\\
	0&0&0&0&w&w&w&w&w&w\\
	0&0&0&0&w&w&w&w&w&w\\ 
	0&0&0&0&w&w&w&w&w&w
\end{pmatrix} \ \hbox{ avec }\  \varphi(1)=\frac{1}{x_1}+O(c),\ \varphi(2)=\frac{1}{x_2}+O(c)\ \hbox{, et }\  \varphi(3)=0.$$
où les  $w$ désignent des expressions algébriques qui ne sont  en général ni nulles ni m\^eme $O(c)$.
\n On en déduit que la connexion tautologique $A$ sur ${\cal E} $ préserve en fait  le sous-fibré vectoriel ${\cal E}'$ de rang 4 de  ${\cal E} $   qui est  engendré par les quatre premiers vecteurs $\bigl(N(j)\bigr)_{ 1\leq j\leq 4}$ de la base $\bigl(N(j)\bigr)_{ 1\leq j\leq 10}$ utilisée pour trivialiser  ${\cal E}$, et que la courbure de la connexion $A'$ ainsi induite sur ${\cal E}'$ est nulle quel que soit $c$. 
Ceci n'a rien d'étonnant :  ${\cal E'}=R'_1$ est en effet le fibré jouant vis-à-vis du 5-sous-tissus $W'_c$ de $W_c$  engendré par les cinq premiers champs de vecteurs le m\^eme r\^ole\footnote{C'est pour qu'il en soit bien ainsi, que -dans la programmation qui suit- ce sont les indices colonne $h_0.i+1+|LL(i-1)|$, $(1\leq i\leq 10)$ que l'on a choisi de supprimer pour obtenir YYY à partir de $M_2$ } que ${\cal E}=R_2$ vis-à-vis de $W_c$, et la restriction $A'$ à ${\cal E'}$ de la connexion tautologique sur ${\cal E}$   est alors la connexion tautologique  sur ${\cal E'}$. Or,  le 5-sous-tissu $W'_c$ est le  5-tissu sur un ouvert de l'espace projectif {\bb P}$^3$ engendré par les 5-bouquets des  droites passant respectivement par les points $[1;0;0;0]$, $[0;1;0;0]$, $[0;0;1;0]$, $[-c;-c;-c;1]$ et $[1;1;1;1]$ : on sait (et l'on redémontre par la m\^eme occasion) qu'un tel tissu a le rang maximum 4 (donc aussi tous les 5-sous-tissus de $W_{0,6}$ puisqu'ils sont tous isomophes), d'où la nullité de la courbure de $A'$.

 Notant de m\^eme $A'(k)$ la sous-matrice $4\times 4$ de $A(k)$ formée avec les 4 premières lignes et les quatre premières colonnes ($A'=\sum_{k=1}^3 A'(k)\ dx[k]$),  on observe que le sous fibré ${\cal E}"$ de rang 1 engendré par $N(1)$ est préservé par la connexion $A'$  sur ${\cal E'}$ ;   plus précisément,   $$\nabla _{\partial_1}N(1)=\Bigl(\frac{1}{x_1}+O(c)\Bigr) N(1) , \hskip 1cm \nabla _{\partial_2}N(1)=\Bigl(\frac{1}{x_2}+O(c)\Bigr)N(1) \hskip 1cm\hbox{ et   }\hskip .5cm
  \nabla _{\partial_3}N(1)=0, $$la connexion $A"$  ainsi induite sur ${\cal E}"$ étant  sans courbure : ceci traduit le fait que le 4-sous-tissu engendré par les $V_i$ et $V_4$ est de rang maximum un  $\Bigl(\bigl({\cal E}"=R_0",A"\bigr)$ jouant, vis à vis du 4-tissu $W_c"$ engendré par les quatre premiers champs de vecteurs (que l'on sait, et que l'on redémontre,  \^etre de rang maximum 1),  le m\^eme r\^ole que 
  $\bigl({\cal E}=R_2,A\bigr)$ vis-à-vis de  $W_c$ ou $\bigl({\cal E}'=R_1,A'\bigr)$ vis-à-vis de  $W_c'\Bigr)$. 

\section {  Appendice (Programmation) :}

$>$ restart;

$>$ with(LinearAlgebra):

$>$ interface(rtablesize=100);

$>$ n:=3;k0:=3;d:=n+k0;

$>$ X:=[seq(x[i],i=1..n)];

$>$  apply(delta,t,s):
for t to n do for s to n do if (t=s) then delta(t,s):=1 else delta(t,s):=0 end if od od;

\n {\it Génération et indexation des multi-indices de dérivation d'ordre 0 à k0+1 }

$>$ apply(L,tau):apply(E,r,y):apply(LL,z):

$>$ for l from 0 to $k0*(k0+1)^{n-2}$ do for k to k0+1 do E(l,k):=Vector(n) od od:

$>$ tau:=1:

$>$ for k to k0+1 do for l from 0 to $k*(k+1)^{n-2}$  do   p:=l:  for s to n-1 do r:=p mod (k+1);

$>$ E(l,k)[s+1]:=r;p:=(p-r)/(k+1) od :  SS:=sum('E(l,k)[u]', 'u'=2..n): 
if SS $<$ (k+1) \hb \indent then E(l,k)[1]:=k-SS:L(tau):=E(l,k):tau:=tau+1 fi od od:

$>$ for i to n do L(0)[i]:=0 od: LL(0):=[seq(L(0)[i],i=1..n)]:

$>$ for t from 0 to binomial(n+k0, k0)-1 do LL(t) := [seq(L(t)[i], i = 1 .. n)] od;

\n {\it fonction  réciproque :}

$>$ apply(ILL, S);  

$>$ for t from 0 to binomial(n+k0, k0)-1 do ILL(LL(t)) := t od;

\n {\it fonction  hauteur 
	d'un multi-indice :}

 $>$ apply(h, t); 

 $>$  for t from 0 to binomial(n+k0, k0)-1 do h(t) := add(L(t)[i], i = 1 .. n) od;

\n  {\it fonctions augmentation 
 d'un multi-indice (dériver une fois de plus par rapport à l'une des coordonnées) :}

$>$ for k to n do for t from 0 to binomial(n+1, k0-1) do apply(ad, k, t) od od;

$>$ apply(ad, k, t); 

$>$ for t from 0 to binomial(n+k0, k0)-1 do for k to n do ad(k, t) := ILL([seq(LL(t)[h]+delta(h, k), h = 1 .. n)]) od od ;

\n {\it Entrée des d champs de vecteurs  (avec introduction d'un  paramètre $c$ de déformation) : } 

 $>$ apply(f,i,a,X):

$>$ for a to n do for i to n do  f(i,a,X):=delta(i,a) od od ;

$>$ for i to n do f(i, n+1, X) := (op(i, X)+c)/(op(n, X)+c) od;

$>$ for i to n do f(i, n+2, X) := (op(i, X)-1)/(op(n, X)-1) od;

$>$ for i to n do f(i, n+3, X) := op(i, X)*(op(i, X)-1)/(op(n, X)*(op(n, X)-1)) od;

{\it  Numérotation des lignes i=(H,u) et colonnes j=(K,v):}

$>$ H :=$ i\to  LL(floor((i+d-2)/(d-1))-1) $; 

$>$ uu :=$ i\to i-(d-1)*(floor((i+d-2)/(d-1))-1) $;

$>$ K := $ j\to LL(floor((j+k0-1)/k0)-1) $; 

$>$ vv := $ i\to j-k0*(floor((j+k0-1)/k0)-1) $;

\n {\it Calcul  de $M_1$ $($H=LL(0), $|K|=0$ ou 1$)$ : }

$>$ apply(F,i,j);

$>$ for i to n do for j to binomial(n+1, 0)*k0 do F(i, j) := simplify(diff(f(i, n+j, X), x[i])) od od;

$>$ for i to n do for j from i*k0+1 to (i+1)*k0 do F(i, j) := f(i, n+vv(j), X); print(

$>$ for i to n do for j from k0+1 to i*k0 do F(i, j) := 0 end do end do; for i to n do for j from (i+1)*k0+1 to (n+1)*k0 do F(i, j) := 0 
od od ;

$>$ for i from n+1 to d-1 do F(i, i-n) := add(diff(f(k, i, X), x[k]), k = 1 .. n) od ;

$>$ for i from n+2 to d-1 do for j to i-n-1 do F(i, j) := 0 od od ; 

$>$ for i from n+1 to d-1 do for j from 0 to binomial(n, 1)-1 do F(i, j*k0+i) := f(j+1, i, X) od od ;

$>$ for i from n+1 to d-1 do for j from i-n+1 to i-n+k0-1 do F(i, j) := 0 od od ;

$>$ for i from n+1 to d-1 do for j from -1 to n-1 do for k from j*k0+i+1 to i+(j+1)*k0-1 do F(i, k) := 0 od od od ;

$>$ MM1 := Matrix(d-1, k0*binomial(n+1, 1), F);

$>$ evalb(Rank(MM1)=d-1);

\n {\it Calcul de $M_2$ $(|H|\leq 1, |K|\leq 2)$ : }

\n  {\it t=0 $($i au plus  $d-1)$} :

$>$ for j from (d-n)*binomial(n+1, 1)+1 to (d-n)*binomial(n+3, 3) do for i to d-1 do F(i, j) := 0 od od ;

\n  {\it $h(t)=1$  :}

$>$ for k to n do for u to d-1 do for v to k0 do F(u+(d-1)*k, v) := simplify(diff(F(u, v), x[k])) od od od;

$>$ for k to n do for u to d-1 do for v to k0 do for s to n do\hb \indent  if s $<>$ k then F(u+(d-1)*k, k0*s+v) := simplify(diff(F(u, k0*s+v), x[k])) \hb \indent  else F(u+(d-1)*k, k*k0+v) := simplify(diff(F(u, k*k0+v), x[k])+F(u, v)) fi od od od od ;

$>$ for k to n do for u to d-1 do for v to k0 do for s from n+1 to binomial(n+3, 3)-1 do if LL(s)[k] = 0 then F(u+(d-1)*k, k0*s+v) := 0  fi od od od od ;

$>$ for k to n do for u to d-1 do for v to k0 do for s from n+1 to binomial(n+3, 3)-1 do if LL(s)[k] $<>$ 0 then F(u+(d-1)*k, k0*s+v) := simplify(F(u, v+k0*sous(k, s))) fi od od od od ;

$>$ MM2 := Matrix((d-1)*binomial((n+1,1), k0*binomial(n+2, 2), F);

$>$ ro :=k0*binomial(n+2, 2)-(d-1)*binomial((n+1,1);

$>$ evalb(Rank(MM2)=(d-1)*binomial((n+1,1));

\n {\it Calcul  de $(Q_3 \ P_3)$ \   $(|H|= 2, |K|\leq 3)$ : }

$>$ for u to d-1 do for v to k0 do for t from binomial(n+1, 1) to binomial(n+2, 2)-1 do 
$F(u+(d-1)*t, v) := simplify(diff(F(u, v), [x_1], [x_2],  LL(t)[2], op(3, X)LL(t)[3]]))$
 end do end do end do;  
 
 \n {\it  $h(t) = 2$ ; pour chaque t, on définit 2 entiers k et m entre 1 et n tels que  $t=ad(k,m)$, avec par exemple $m\leq k$}
 
 $>$ apply(phi,t):apply(psi,t):
 
 $>$ phi(4):=1 :phi(5):=2 :phi(6):=2 :phi(7):=3 :phi(8):=3 :phi(9):=3 :
 
 $>$ psi(4):=1 :psi(5):=1 :psi(6):=2 :psi(7):=1 :psi(8):=2 :psi(9):=3 :

\n {\it Vérification :}

$>$ for t from n+1 to binomial(n+2,2)-1 do evalb(ad(phi(t),psi(t)=t)) od;

\n {\it si $L(s)[phi(t)]=0$ :}

 $>$ for t from n+1 to binomial(n+2,2)-1 do for s from 0 to binomial(n+3,3)-1 do\hb 
 if L(s)[phi(t)]=0 then for u to d-1 do for v to k0 do F(u+(d-1)*t,v+k0*s):=simplify(diff(F(u+(d-1)*psi(t),v+k0*s),x[phi(t)])) od od fi od od :

\n {\it si s de la forme $s=ad(phi(t),r)$ :}

$>$ for t from n+1 to binomial(n+2,2)-1 do for r from 0 to binomial(n+2,2)-1 do 
 for u to d-1 do for v to k0 do F(u+(d-1)*t,v+k0*ad(phi(t),r):=simplify(diff(F(u+(d-1)*psi(t),v+k0*ad(phi(t),r)),x[phi(t)])+F(u+(d-1)*psi(t),v+k0*r) od od  od od :

\hskip .3cm

$>$ apply(qq, i, j); for i to 30 do for j to 30 do qq(i, j) := F(i+20, j) od od;

$>$Q3 := Matrix((d-1)*binomial(n+1,2), k0*binomial(n+2, 2), F);

$>$ apply(p3, i, j); for i to 30 do for j to 30 do p3(i, j) := F(i+20, j+30) od od;

$>$ P3 := Matrix((d-1)*binomial(n+1,2), k0*binomial(n+2, 3), p3);

$>$ evalb(Rank(P3)=(d-1)*binomial(n+1,2));

$>$ IP3 := simplify(MatrixInverse(P3));

$>$ U := simplify(-IP3 . Q3);

\n {\it Définition d'une sous-matrice carrée inversible YYY de $M_2$:}

{\it liste  des indices colonne de $M_2$ omis pour obtenir YYY:}

$>$ for i to ro do $k0*(i-1)+1+h(i-1)$  od;

\hskip 5cm 1 5 8 11 15 18 21 24 27 30

{\it Ces indices sont choisis, de façon que les quatre premiers d'entre eux (1,5,8,11) soient ceux que l'on obtiendrait par la m\^eme procédure avec le 5-sous-tissu obtenu en supprimant le sixième feuilletage (cf. remarques finales précédant l'appendice).}

 {\it - liste  complémentaire des indices colonnes   conservés }
 
 $>$ for j to (d-1)*binomial(n+1,1) do $floor(\frac{3*j+1-h(floor(\frac{j-1}{2}))}{2})$ od ;
 
\hskip 3cm 2 3 4 6 7 9 10 12 13 14 16 17 19 20 22 23 25 26 28 29

$>$ apply(yyy,i,j);
 
$>$  for i to (d-1)*binomial(n+1,1) do for j to (d-1)*binomial(n+1,1)   do $yyy(i,j):=F[i,floor(\frac{3*j+1-h(floor(\frac{j-1}{2}))}{2})]$ od od;

$>$  YYY:=Matrix((d-1)*binomial(n+1,1),(d-1)*binomial(n+1,1),yyy)

$>$ evalb(Rank(YYY)=(d-1)*binomial(n+1,1));

\n {\it Calcul d'une base N(j) , j=1..10, de l'espace des sections de $\cal E$ := Ker $(M_2)$ :}

$>$ W := Vector(k0*binomial(n+2,2), symbol = w);

 $>$ EE := simplify(MM2 . W);
 
$>$ S := simplify(solve({EE[1], EE[2], EE[3], EE[4], EE[5], EE[6], EE[7], EE[8], EE[9], EE[10], EE[11], EE[12], EE[13], EE[14], EE[15], EE[16], EE[17], EE[18], EE[19], EE[20]}, {w[2], w[3], w[4],  w[6], w[7],  w[9], w[10], w[12], w[13],w[14],w[16],w[17],w[19],w[20], w[22], w[23],  w[25], w[26], w[28], w[29]}));

$>$ S := simplify(S);

$>$ assign(S);

$>$ apply(N, j);

$>$ N(1) := simplify(subs({w[1] = 1, w[5] = 0, w[8] = 0, w[11] = 0, w[15] = 0, w[18] = 0, w[21] = 0, w[24] = 0, w[27] = 0, w[30] = 0}, [w[1], w[2], w[3], w[4], w[5], w[6], w[7], w[8], w[9], w[10], w[11], w[12], w[13], w[14], w[15], w[16], w[17], w[18], w[19], w[20], w[21], w[22], w[23], w[24], w[25], w[26], w[27], w[28], w[29], w[30]]));

$>$ N(2) := simplify(subs({w[1] = 0, w[5] = 1, w[8] = 0, w[11] = 0, w[15] = 0, w[18] = 0, w[21] = 0, w[24] = 0, w[27] = 0, w[30] = 0}, [w[1], w[2], w[3], w[4], w[5], w[6], w[7], w[8], w[9], w[10], w[11], w[12], w[13], w[14], w[15], w[16], w[17], w[18], w[19], w[20], w[21], w[22], w[23], w[24], w[25], w[26], w[27], w[28], w[29], w[30]]));

$>$ N(3) := simplify(subs({w[1] = 0, w[5] = 0, w[8] = 1, w[11] = 0, w[15] = 0, w[18] = 0, w[21] = 0, w[24] = 0, w[27] = 0, w[30] = 0}, [w[1], w[2], w[3], w[4], w[5], w[6], w[7], w[8], w[9], w[10], w[11], w[12], w[13], w[14], w[15], w[16], w[17], w[18], w[19], w[20], w[21], w[22], w[23], w[24], w[25], w[26], w[27], w[28], w[29], w[30]]));

$>$ N(4) := simplify(subs({w[1] = 0, w[5] = 0, w[8] = 0, w[11] = 1, w[15] = 0, w[18] = 0, w[21] = 0, w[24] = 0, w[27] = 0, w[30] = 0}, [w[1], w[2], w[3], w[4], w[5], w[6], w[7], w[8], w[9], w[10], w[11], w[12], w[13], w[14], w[15], w[16], w[17], w[18], w[19], w[20], w[21], w[22], w[23], w[24], w[25], w[26], w[27], w[28], w[29], w[30]]));

$>$ N(5) := simplify(subs({w[1] = 0, w[5] = 0, w[8] = 0, w[11] = 0, w[15] = 1, w[18] = 0, w[21] = 0, w[24] = 0, w[27] = 0, w[30] = 0}, [w[1], w[2], w[3], w[4], w[5], w[6], w[7], w[8], w[9], w[10], w[11], w[12], w[13], w[14], w[15], w[16], w[17], w[18], w[19], w[20], w[21], w[22], w[23], w[24], w[25], w[26], w[27], w[28], w[29], w[30]]));

$>$ N(6) := simplify(subs({w[1] = 0, w[5] = 0, w[8] = 0, w[11] = 0, w[15] = 0, w[18] = 1, w[21] = 0, w[24] = 0, w[27] = 0, w[30] = 0}, [w[1], w[2], w[3], w[4], w[5], w[6], w[7], w[8], w[9], w[10], w[11], w[12], w[13], w[14], w[15], w[16], w[17], w[18], w[19], w[20], w[21], w[22], w[23], w[24], w[25], w[26], w[27], w[28], w[29], w[30]]));

$>$ N(7) := simplify(subs({w[1] = 0, w[5] = 0, w[8] = 0, w[11] = 0, w[15] = 0, w[18] = 0, w[21] = 1, w[24] = 0, w[27] = 0, w[30] = 0}, [w[1], w[2], w[3], w[4], w[5], w[6], w[7], w[8], w[9], w[10], w[11], w[12], w[13], w[14], w[15], w[16], w[17], w[18], w[19], w[20], w[21], w[22], w[23], w[24], w[25], w[26], w[27], w[28], w[29], w[30]]));

$>$ N(8) := simplify(subs({w[1] = 0, w[5] = 0, w[8] = 0, w[11] = 0, w[15] = 0, w[18] = 0, w[21] = 0, w[24] = 1, w[27] = 0, w[30] = 0}, [w[1], w[2], w[3], w[4], w[5], w[6], w[7], w[8], w[9], w[10], w[11], w[12], w[13], w[14], w[15], w[16], w[17], w[18], w[19], w[20], w[21], w[22], w[23], w[24], w[25], w[26], w[27], w[28], w[29], w[30]]));

$>$ N(9) := simplify(subs({w[1] =0, w[5] = 0, w[8] = 0, w[11] = 0, w[15] = 0, w[18] = 0, w[21] = 0, w[24] = 0, w[27] = 1, w[30] = 0}, [w[1], w[2], w[3], w[4], w[5], w[6], w[7], w[8], w[9], w[10], w[11], w[12], w[13], w[14], w[15], w[16], w[17], w[18], w[19], w[20], w[21], w[22], w[23], w[24], w[25], w[26], w[27], w[28], w[29], w[30]]));

$>$ N(10) := simplify(subs({w[1] = 0, w[5] = 0, w[8] = 0, w[11] = 0, w[15] = 0, w[18] = 0, w[21] = 0, w[24] = 0, w[27] = 0, w[30] = 1}, [w[1], w[2], w[3], w[4], w[5], w[6], w[7], w[8], w[9], w[10], w[11], w[12], w[13], w[14], w[15], w[16], w[17], w[18], w[19], w[20], w[21], w[22], w[23], w[24], w[25], w[26], w[27], w[28], w[29], w[30]]));

$>$ apply(n2, i, j); $>$ for i to 30 do for j to 10 do n2(i, j) := N(j)[i] end do end do;

$>$ N2 := Matrix(k0*binomial(n+2, 2), r0, n2);

\n {\it Vérification : MM2.N2=0 et rang(N2)=10 }

$>$ simplify(MM2 . N2);

$>$ Rank(N2);

\n {\it Dérivées covariantes :}

$>$ n3:=simplify(U.N2):

$>$ for k to n do for s from 0 to n do for v to k0 do for j to ro do \hb nn3(k)(k0*s+v, j) := N2[v+k0*ad(k, s), j] od od od od ; 

$>$ for k to n do for s from n+1 to binomial(n+2, 2)-1 do for v to k0 do for j to ro do \hb nn3(k)(k0*s+v, j) := n3[v+k0*ad(k, s)-k0*binomial(n+2, 2), j] od od od od;

$>$ apply(N3, k); for k to n do N3(k) := Matrix(k0*binomial(n+2, 2), ro, nn3(k)) od;

$>$ apply(ff, k);

$>$  for k to n do apply(ff(k), i, j) od ;

$>$ for k to n do for j to ro do for i to (d-n)*binomial(n+2, 2) do (ff(k))(i, j) := simplify(diff(N(j)[i], x[k])) od od od ;

$>$ apply(dN2, k);

$>$ for k to n do dN2(k) := Matrix((d-n)*binomial(n+2, 2), ro, ff(k)) od;

$>$ (dc, k); apply(DC, k);

$>$ for k to n do apply(dc(k), i, j) od;

$>$ for k to n do for i to (d-n)*binomial(n+2, 2) do for j to ro do\hb  dc(k) :=(i, j) $\to$ simplify((ff(k))(i, j)-(nn3(k))(i, j))  od od od ;

$>$ for k to n do DC(k) := Matrix((d-n)*binomial(n+2, 2), ro, dc(k)) od;

\n {\it Formes de connexion :
on ne garde que  les lignes  de DC(k) dont les indices sont ceux des colonnes omises pour obtenir YYY à partir de MM2 :}

$>$ apply(A, k); apply(alpha, k); for k to n do apply(alpha(k), i, j) od;

$>$ for k to n do for i to (d-n)*binomial(n+2, 2) do for j to ro do \hb alpha(k)(i, j) := simplify((dc(k))(k0*(i-1)+1+h(i-1), j)) od od od ;

$>$ for k to n do A(k) := Matrix(ro, ro, alpha(k)) od;

\n {\it Courbure :}

$>$ apply(dalpha, k, m); apply(b, k, m); apply(e, k, m); apply(braket, k, m);

$>$ for k to n-1 do for m from k+1 to n do apply(dalpha(k, m), i, j) od od ;

$>$ for k to n-1 do for m from k+1 to n do for i to ro do for j to ro do \hb
dalpha(k, m)(i, j) := simplify(simplify(diff(A(m)[i, j], x[k]))-simplify(diff(A(k)[i, j], x[m]))) od od od od ;

$>$  for k to n-1 do for m from k+1 to n do apply(b(k, m), i, j) od od ; 

$>$ for k to n-1 do for m from k+1 to n do for i to ro do for j to ro do \hb b(k, m)(i, j) := simplify((A(k) . A(m))[i, j]) od od od od ;

$>$ for k to n-1 do for m from k+1 to n do apply(e(k, m), i, j) od od ; 

$>$ for k to n-1 do for m from k+1 to n do for i to ro do for j to ro do \hb e(k, m)(i, j) := simplify((A(m) . A(k))[i, j]) od od od od ;

$>$ for k to n-1 do for m from k+1 to n do apply(braket(k, m), i, j) od od ;

$>$ for k to n-1 do for m from k+1 to n do for i to ro do for j to ro do \hb
braket(k, m)(i, j) := simplify(simplify((b(k, m))(i, j))-simplify((e(k, m))(i, j))) od od od od ;

$>$ apply(ko, k, m); for k to n-1 do for m from k+1 to n do apply(ko(k, m), i, j) od od ;

$>$ for k to n-1 do for m from k+1 to n do for i to ro do for j to ro do \hb 
ko(k, m)(i, j) := simplify(dalpha(k, m)(i, j)+braket(k, m)(i, j)) od od od od ;

$>$ for k to n-1 do for m from k+1 to n do for i to ro do for j to ro do \hb koO(k, m)(i, j) := taylor(ko(k, m)(i, j), c, 1) od od od od ;

$>$ apply(KoO, k, m);

$>$ for k to n-1 do for m from k+1 to n do KoO(k, m) := Matrix(ro, ro, koO(k, m)); print(Ko(k, m) =KoO(k,m)) od od ;

\n {\it On vérifie que la connexion préserve le sous-fibré engendré par $\bigl(N(1),..., N(4)\bigr)$ : tous les coefficients des matrices $\Psi(k)$   doivent \^etre nuls.}

$>$ for k to n do beta(k):=(i,j)$\to$ taylor(alpha(k)(i,j),c,1) od;  

$>$ for k to n do A(k):=Matrix(ro,ro,beta(k)) od;

$>$ for k to n do $\Psi(k)$=DeleteRow((DeleteColumn(A(k),5..10)),1..4) od;

$>$ for k to n do $A'(k)$=DeleteRow((DeleteColumn(A(k),1..4)),1..4) od;
\hskip 2mm
 
\n {\bf Références}
	
\hskip 2mm
	
\n [BB] W. Blaschke et G. Bol, {\it Geometrie der Gewebe}, Die Grundlehren der Mathematik 49, Springer, 1938.
\hskip 1mm
	
\n [Bo] G. Bol, {\it \"Uber ein bemerkenswertes F\"unfgewebe in der Ebene}, Abh. Math. Hamburg Univ., 11, 1936, 387-393.
\hskip 1mm
	
	
	
\n [CL] V. Cavalier, D. Lehmann, {\it Ordinary holomorphic webs of codimension one. } arXiv 0703596v2 [mathsDS], 2007, et Ann. Sc. Norm. Super. Pisa, cl. Sci (5), vol XI (2012), 197-214.
\hskip 1mm
	
\n [D1] D.B. Damiano, {\it Abelian equations and characteristic classes}, Thesis, Brown University, (1980) ;  American J. Math. 105-6, 
	1983,  1325-1345.
	\hskip 1mm

\n [D2] D.Damiano : {\it Webs  and characteristic forms on Grassmann manifolds},Am.J. of Maths.105, 1983,  1325-1345.
\hskip 1mm

\n [DL1] J. P. Dufour, D. Lehmann, {\it Calcul explicite de la courbure des tissus calibr\'es ordinaires,; } arXiv 1408.3909v1 [mathsDG], 18/08/2014.
	\hskip 1mm
	
\n [DL2] J. P. Dufour, D. Lehmann, {\it Rank of ordinary webs in codimension one : an effective method ; }\hb  arXiv 1703.03725v1 [math.DG], 10/03/2017. Pure and Applied Mathematics Quarterly, vol. 16, n. 5, 1587-1607, 2020. 
	\hskip 1mm
	
\n [DL3] J. P. Dufour, D. Lehmann, {\it Rank of ordinary webs in codimension one : an effective method}. Pure and Applied Mathematics Quarterly, vol. 16-n°5 (2020), 1587-1607.

	\n [DL4] J.P. Dufour, D. Lehmann, {\it Etude des $(n+1)$-tissus de courbes en dimension $n$}, arXiv 2211.05188v1 [mathsDG], 09/11/2022, et 	Comptes Rendus Maths.  Ac. Sc. Paris, vol. 361, 1491-1497, 2023.
	\hskip 1mm 
	
	

\n [H1] A. Hénaut, {\it Planar web geometry through abelian relations  and connections}
	Annals of Math. 159 (2004)  425-445.
		\hskip 1mm
	
\n [H2]  A. Hénaut, {\it Formes différentielles abéliennes, bornes de Castelnuovo et géométrie des tissus}, Commentarii Math.  Helvetici, 79 (1), 2004, 25-57.
	\hskip 1mm
	
\n [L]  D. Lehmann, {\it Relations abéliennes des tissus ordinaires de codimension arbitraire}, arXiv:1712.00997,\hb v1(4/12/2017),v2(30/12/2021).
\hskip 1mm

\n [Pa] A. Pantazi. {\it Sur la d\'etermination du rang d'un tissu plan.} C.R. Acad. Sc. Roumanie 4 (1940), 108-111.
	\hskip 1mm

\n [Pi1] L. Pirio, {\it Equations Fonctionnelles Ab\'eliennes et G\'eom\'etrie des tissus},
	Th\`ese de doctorat de l'Universit\'e Paris VI, 2004.
		\hskip 1mm

\n [Pi2]  L. Pirio, {\it Sur les tissus planaires de rang maximal et le problème de Chern}, note aux C.R. Ac Sc. , sér. I, 339 (2004), 131-136.
	\hskip 1mm


\n  [Pi3] L. Pirio : {\it On the $(n+3)$-webs  by rational curves induced by the forgetful maps on the moduli spaces ${\cal
			M}_{0,n+3}$}, arXiv 2204.04772.v1, [Math AG], 10-04-2022.
			\vskip 1cm

	\noindent Daniel Lehmann, ancien  professeur \`a l'Universit\'e de {\cal M}ontpellier II, \hfill\break   4 rue Becagrun,  30980 Saint Dionisy, France \hfill\break  email : lehm.dan@gmail.com,

\end{document}